\documentclass[12pt]{article}
\usepackage[a4paper,left=4cm,right=1.2cm,top=3cm,bottom=4cm]
{geometry}
\usepackage{amsmath,amstext, amsthm, amsbsy,amssymb,marvosym,fancyhdr,amscd,amsfonts,latexsym,delarray,stackrel,lineno,color,cite}
\usepackage{graphicx}

\usepackage{wasysym}
\usepackage{srcltx}
\usepackage{mathrsfs}

\usepackage[all]{xy}
\usepackage{t1enc}
\usepackage{mathrsfs}
\usepackage{pifont}

\usepackage{srcltx}
\usepackage{mathrsfs}
\usepackage{mathptmx}       
\usepackage{helvet}         
\usepackage{courier}        
\usepackage{type1cm}        
\usepackage[bottom]{footmisc}
\usepackage{float}  

\newcommand*\patchAmsMathEnvironmentForLineno[1]{%
  \expandafter\let\csname old#1\expandafter\endcsname\csname #1\endcsname
  \expandafter\let\csname oldend#1\expandafter\endcsname\csname end#1\endcsname
  \renewenvironment{#1}%
     {\linenomath\csname old#1\endcsname}%
     {\csname oldend#1\endcsname\endlinenomath}}%
\newcommand*\patchBothAmsMathEnvironmentsForLineno[1]{%
  \patchAmsMathEnvironmentForLineno{#1}%
  \patchAmsMathEnvironmentForLineno{#1*}}%
\AtBeginDocument{%
\patchBothAmsMathEnvironmentsForLineno{equation}%
\patchBothAmsMathEnvironmentsForLineno{align}%
\patchBothAmsMathEnvironmentsForLineno{flalign}%
\patchBothAmsMathEnvironmentsForLineno{alignat}%
\patchBothAmsMathEnvironmentsForLineno{gather}%
\patchBothAmsMathEnvironmentsForLineno{multline}%
}

\definecolor{Green}{rgb}{0,1,0}
\definecolor{Blue}{RGB}{0,0,191}
\definecolor{mathmodecolor}{RGB}{0,102,0}
\definecolor{keywordcolor}{RGB}{0,51,151}
\definecolor{sourcebackgroundcolor}{RGB}{255,247,223}
\definecolor{unixagred}{RGB}{255,0,0}
\definecolor{green}{RGB}{1,191,191}
\definecolor{lightgray}{gray}{0.88}

\newtheorem{thm}{Theorem}[section]

\newtheorem{lem}[thm]{Lemma}

\newtheorem{defn}[thm]{Definition}

\def\qqq{\,,\quad~\forall}

\def\Aut{{\rm Aut}}

\def\GL{{\rm GL}}

\def\Hom{{\rm Hom}}
\def\id{{\rm id}}

\def\Ker{{\rm Ker}}

\def\mmod{{\rm Mod}}

\def\Spec{{\rm Spec\,}}

\def\Tr{{\rm Tr}}

\def\A{{\mathbb A}}
\def\B{{\mathbb B}}
\def\C{{\mathbb C}}
\def\F{{\mathbb F}}
\def\K{{\mathbb K}}
\def\N{{\mathbb N}}
\def\P{{\mathbb P}}
\def\Q{{\mathbb Q}}
\def\R{{\mathbb R}}
\def\Z{{\mathbb Z}}

\def\gam{{\Gamma\Ses}}

\newcommand{\har}{HC^{\rm ar}}
\newcommand{\hes}{TP}

\def\urep{\vartheta}

\def\aarith{{\mathscr A}}
\def\scal{{(\rnt,\cO)}}
\def\scal1{{\hat \aarith}}
\def\hatz{{\hat\Z^\times}}

\def\mapE{{\mathfrak E}}
\def\ett{{\rm et}}
\def\Tr{{\rm Tr}}

\def\cdim{{{\mbox{Dim}_\R}}}
\def\tdim{{{\mbox{dim}_{\rm top}}}}
\def\inter{{\mathfrak s}}
\def\cA{{\mathcal A}}

\def\cE{{\mathcal E}}

\def\cH{{\mathcal H}}

\def\cK{{\mathcal K}}

\def\cM{{\mathcal M}}
\def\cO{{\mathcal O}}
\def\cP{{\mathcal P}}

\def\cS{{\mathcal S}}

\def\qqq{\,,\,~\forall}

\newcommand{\ie}{{\it i.e.\/}\ }
\newcommand{\eg}{{\it e.g.\/}\ }
\newcommand{\cf}{{\it cf.}}
\newcommand{\cff}{{\it cf.}~}
\newcommand{\opcit}{{\it op.cit.\/}\ }

\def\spz{{\Spec\Z}}

\def\zmin{\Z_{\rm min}}

\def\id{{\mbox{Id}}}

\def\dim{{\mbox{dim}}}

\def\Hom {{\mbox{Hom}}}

\def\catmo{{\mathfrak{Mod}}}

\def\mc{multiplicatively cancellative }

\def\rma{\R_{\rm max}}
\def\sss{{\mathbb  S}}
\def\gop{{\Gamma^{\rm op}}}

\def\zmax{{\Z_{\rm max}}}

\def\rmax{\R_+^{\rm max}}

\def\fr{{\rm Fr}}

\def\arith{{(\wnt,\bar \N)}}

\def\Deg{{\rm Deg}}

\def\vsp{\vspace{.05in}}

\def\Se{\frak{ Sets}}

\newcommand{\nil}[1]{}

\def\dop{{\Delta^{\rm op}}}

\def\nt{\N^{\times}}
\def\wnt{{\widehat{\N^{\times}}}}
\def\rnt{{[0,\infty)\rtimes{\N^{\times}}}}

\def\div{{\rm Div}}

\def\nbo{{\zmin\otimes_\B \zmin}}

\def\arith{{(\wnt,\zmax)}}
\def\wntb{{\widehat{\N^{\times 2}}}}

\def\Ses{{\Se_*}}

\def\li{{\rm Li}}

\def\rmz{{\rm Zeros}}

\begin{document}

\title{An essay on the Riemann Hypothesis}
\author{Alain Connes\footnote{Coll\`ege de France, 3 rue d'Ulm, Paris 75005 France. 
IH\'ES, 35 Route de Chartres, Bures sur Yvette  and Ohio State University, Columbus, Ohio; email : alain@connes.org}}

%
%

\maketitle

\abstract{The Riemann hypothesis is, and will hopefully remain for a long time, a great motivation to uncover and explore new parts of the mathematical world. After reviewing its impact on the development of algebraic geometry we discuss three strategies, working concretely at the level of the explicit formulas. The first strategy is ``analytic" and is based on Riemannian spaces and Selberg's work on the trace formula and its comparison with the explicit formulas. The second is based on algebraic geometry and the Riemann-Roch theorem. We establish a framework 
in which one can transpose many of the ingredients of the Weil proof as reformulated by Mattuck, Tate and Grothendieck. This framework is elaborate and involves noncommutative geometry, Grothendieck toposes and tropical geometry. We point out the remaining difficulties and show that RH gives a strong motivation to develop algebraic geometry in the emerging world of characteristic one. Finally we briefly discuss a third strategy based on the development of a suitable ``Weil cohomology",  the role of Segal's $\Gamma$-rings and of topological cyclic homology as a model for ``absolute algebra" and as a cohomological tool.\newline\indent
}

\tableofcontents
\section{Introduction}
\label{intro}
Let $\pi(x):=\#\{p\mid p\in \cP, \, p<x\}$ be the number of primes less than $x$ with $\frac12$
added when $x$ is prime.
Riemann \cite{Riemann} found  for the counting function \footnote{Similar counting functions were already present in Chebyshev's work}
$$
f(x):=\sum \frac 1n \pi(x^{\frac 1n}), 
$$
the following formula
 involving the integral logarithm function $\li(x)=\int_0^x\frac{dt}{\log t}$,
\begin{equation}\label{Riemann1}
f(x)=\li(x)-\sum_\rho \li(x^{\rho})+\int_x^\infty \frac{1}{t^2-1}\,\frac{dt}{t\log t}-\log 2
\end{equation}
in terms\footnote{More precisely Riemann  writes $\sum_{\Re(\alpha)>0}\left(\li(x^{\frac 12+\alpha i})+\li(x^{\frac 12-\alpha i})\right)$ instead of $\sum_\rho \li(x^{\rho})$ using the symmetry $\rho\to 1-\rho$
provided by the functional equation, to perform the summation.} of the non-trivial zeros  $\rho$ of the analytic continuation (shown as well as two proofs of the functional equation by Riemann at the beginning of his paper) of the Euler zeta function
$$
\zeta(s)=\sum\frac{1}{n^s}
$$
Reading Riemann's original paper is surely still the best initiation to the subject.    In  his lecture given in Seattle  in August 1996, on the occasion of the 100-th anniversary of the proof of the prime number theorem, Atle Selberg comments about Riemann's paper:           \cite{Selbergseattle}
\begin{quote} {\em It is clearly a preliminary note and might not have been written if L. Kronecker had not urged him to write up something about this work (letter to Weierstrass, Oct. 26 1859). It is clear that there are holes that need to be filled in, but also clear that he had a lot more material than what is in the note\footnote{See \cite{RZeta} Chapter VII for detailed support to Selberg's comment}. What also seems clear : Riemann is not interested in an asymptotic formula, not in the prime number theorem, what he is after is an exact formula!}
\end{quote}

The Riemann hypothesis (RH) states that all the non-trivial zeros of $\zeta$ are on the line $\frac 12+i\R$. This hypothesis has become over the years and the many unsuccessful attempts at proving it, a kind of ``Holy Grail" of mathematics. Its validity is indeed one of the deepest conjectures and besides its clear inference on the distribution of prime numbers, it admits relations with many parts of pure mathematics as well as of quantum physics.

 It is, and will hopefully remain for a long time, a great motivation to uncover and explore new parts of the mathematical world. 
There are many excellent texts on RH, such as \cite{B2} which explain in great detail what is known about the problem, and the many implications of a positive answer to the conjecture. When asked by John Nash to write a text on RH\footnote{My warmest thanks to Michael Th. Rassias for the communication}, I realized that writing  one more encyclopedic text would just add another layer to the psychological barrier that surrounds  RH.  Thus I have chosen deliberately to adopt another point of view, which is to navigate between the many forms of the explicit formulas (of which \eqref{Riemann1} is the prime example) and possible strategies to attack the problem, stressing the value of the elaboration of new concepts rather than  ``problem solving". 

\begin{itemize}
\item {\it RH and algebraic geometry}

\vspace{.05in}
We first explain the Riemann-Weil explicit formulas in the framework of adeles and global fields in \S \ref{subsecrw}. We then sketch in \S \ref{rrstrat} the geometric proof of RH for function fields as done by Weil, Mattuck, Tate and Grothendieck. We then turn to the role of 
RH in generating new mathematics, its role in the evolution of algebraic geometry in the XX-th century through the Weil conjectures, proved by Deligne, and the elaboration by Grothendieck of the notions of scheme and of topos.

\item{\it Riemannian Geometry, Spectra  and trace formulas} 

\vsp
Besides the proof of analogues of RH such as the results of Weil and of Deligne, there is another family of results that come pretty close. They give another natural approach of RH using analysis, based on the pioneering work of Selberg on trace formulas. These will be reviewed in Section \ref{analysisattack} where the difficulty arising from the minus sign in front of the oscillatory terms will be addressed.

\item{\it The Riemann-Roch strategy: A Geometric Framework} 

\vsp 
In Section \ref{algeomattack}, we shall describe a geometric framework, established in our joint work with C. Consani, allowing us to transpose several of the key ingredients of the geometric proof of RH for function fields recalled in \S \ref{rrstrat}. It is yet unclear if this is the right set-up for the final Riemann-Roch step, but  it will illustrate  the power of RH as an incentive to explore new parts of mathematics since it gives a clear motivation for developing algebraic geometry in characteristic $1$ along the line of tropical geometry. This will take us from the world of characteristic $p$ to the world of characteristic $1$, and give us an opportunity to describe its relation with semi-classical and idempotent analysis,  optimization and game theory\footnote{one of the topics in which John Nash made fundamental contributions}, through the Riemann-Roch theorem in tropical geometry \cite{BN,GK,MZ}.

\item {\it Absolute Algebra and the sphere spectrum} 

\vsp 
 The arithmetic and scaling sites which are the geometric spaces underlying the Riemann-Roch strategy of Section \ref{algeomattack} are only the semiclassical shadows of a more mysterious structure underlying the compactification of $\Spec\Z$ that should give a cohomological interpretation of the explicit formulas.  We describe in this last section an essential tool  coming from algebraic topology: Segal's $\Gamma$-rings and the sphere spectrum, over which all previous attempts at developing an absolute algebra organize themselves.  Moreover, thanks to the results of Hesselholt and Madsen in particular, topological cyclic homology gives a cohomology theory suitable to treat in a unified manner the local factors of $L$-functions.

\end{itemize}

\section{RH and algebraic geometry}
I will briefly sketch here the way RH, once transposed in finite characteristic, has played a determining role in the upheaval of the very notion of geometric space in algebraic geometry culminating with the notions of scheme and topos due to Grothendieck, with the notion of topos offering a frame of thoughts of incomparable generality and breadth. It is a  quite remarkable testimony to the unity of mathematics that the origin of this discovery lies in the greatest problem of analysis and arithmetic.

\subsection{The Riemann-Weil explicit formulas, Adeles and global fields}
\label{subsecrw}
Riemann's formula \eqref{Riemann1} is a special case of the ``explicit formulas" which establish a duality between the primes and the zeros of zeta. This formula has been extended by Weil in the context of global fields which provides a perfect framework for a generalization of RH since it has been solved, by Weil, for all global fields except number fields. 

\subsubsection{The case of $\zeta$}

Let us start with the explicit formulas (\cff \cite{weilpos0,weilpos, EB, patter}). We start with a
function
$F(u)$ defined for $u\in [1,\infty)$, continuous and continuously differentiable except for finitely many points at which both $F(u)$ and $F'(u)$ have at most a discontinuity of the first kind, \footnote{and at which the value of $F(u)$ is defined as the average of the right and left limits there} and such that, for some $\epsilon>0$, $F(u)=O(u^{-1/2-\epsilon})$. One then defines the Mellin transform of $F$ as 
\begin{equation}\label{expl}
\Phi(s)=\int_1^\infty F(u)\,u^{s-1}du
\end{equation}
The explicit formula then takes the form
\begin{equation}\label{explfor}
\Phi(\frac 12)+\Phi(-\frac 12)-\sum_{\rho\in \rmz}\Phi(\rho-\frac 12)=\sum_p\sum_{m=1}^\infty \log p \,\,p^{-m/2}F(p^m)+
\end{equation}
$$
+(\frac \gamma 2+\frac{\log\pi}{2})F(1)
+\int_1^\infty\frac{t^{3/2}F(t)-F(1)}{t(t^2-1)}dt
$$
where $\gamma=-\Gamma'(1)$ is the Euler constant, and  the zeros are counted with their multiplicities \ie $\sum_{\rho \in\rmz}\Phi(\rho-\frac 12)$ means $\sum_{\rho\in \rmz}{\rm order}(\rho)\Phi(\rho-\frac 12)$.

\subsubsection{Adeles and global fields}
By a result of Iwasawa \cite{Iwasawa} a field $\K$ is a finite algebraic number field, or an algebraic function field of one variable over a finite constant field, if and only if there exists 
 a semi-simple (\ie with trivial Jacobson radical \cite{Jac}) commutative ring $R$ containing $\K$ such that $R$ is locally compact, but neither compact nor discrete and  $\K$ is discrete and cocompact in $R$. This result gives a conceptual definition of what is a ``global field" and indicates that the arithmetic of such fields is intimately related to analysis on the parent ring $R$ which is called the ring of adeles of $\K$ \cite{Weil,tate}. It is the opening door to a whole world which is that of automorphic forms and representations, starting in the case of $\GL_1$ with Tate's thesis \cite{tate} and Weil's book \cite{Weil}. Given a global field $\K$, the ring  $\A_\K$ of adeles of $\K$  is the restricted product of the locally compact fields $\K_v$ obtained as completions of $\K$ for the different places $v$ of $\K$. The equality $dax=\vert a\vert dx$ for the additive Haar measure defines the module $\mmod:\K_v\to \R_+$, $\mmod(a):=\vert a\vert$  on the local fields $\K_v$ and also as a group homomorphism $\mmod:C_\K\to \R_+^*$ where $C_\K=\GL_1(\A_\K)/\K^\times$ is the idele class group. The kernel of the module is a compact subgroup $C_{\K,1}\subset C_\K$ and the range of the module is  a cocompact subgroup $\mmod(\K)\subset  \R_+^*$. On any locally compact modulated group, such as $C_\K$ or the multiplicative groups $\K_v^*$, one normalizes the Haar measure $d^*u$ uniquely so that the measure of $\{u\mid 1\leq \vert u\vert \leq \Lambda\}$ is equivalent to $\log \Lambda$ when $\Lambda\to \infty$.

\subsubsection{Weil's explicit formulas}\label{sectweilexpl}

As shown by Weil, in \cite{weilpos}, adeles and global fields give the natural framework for the explicit formulas.  For each character $\chi\in \widehat{C_{\K,1}}$ one chooses an extension $\tilde \chi$ to $C_\K$ and one lets $Z_{\tilde\chi}$ be the set (with multiplicities and taken modulo the orthogonal of $\mmod(\K)$, \ie $\{s\in\C\mid q^s=1, \forall q\in \mmod(\K)\}$) of zeros of the $L$-function associated to $\tilde\chi$.
Let  then $ \alpha$ be a nontrivial character of 
 $\A_\K/\K$ and $ \alpha = \prod \,  \alpha_v$ its local factors. The explicit formulas take the following form, with $h \in \cS (C_\K)$ a Schwartz function with compact support:
\begin{equation}\label{weil4}
\hat h (0) + \hat h (1)  - \sum_{\chi\in \widehat{C_{\K,1}}}\,
\sum_{Z_{\tilde\chi}} \hat h (\tilde\chi , \rho) = \sum_v
\int'_{\K_v^*} \frac{h(u^{-1}) }{   \vert 1-u   \vert} \, d^* u
\end{equation}
where  the principal value $\int'_{\K_v^*}$ is normalized by the additive character $ \alpha_v$ (\cff\cite{CMbook} Chapter II, 8.5, Theorem 2.44 for the precise notations and normalizations) and for any character $\omega$ of $C_\K$ one lets
\begin{equation}\label{fourier}
\hat h (\omega , z): = \int h(u) \,
\omega(u) \,   \vert u   \vert^z \, d^* u, \  \  \hat h (t):=\hat h (1,t)
\end{equation}
For later use in \S \ref{counting} we compare \eqref{explfor} with the Weil way \eqref{weil4} of writing the explicit formulas. Let the function $h$ be the function on $C_\Q$ given by $h(u):=\vert u\vert^{-\frac 12}F(\vert u\vert)$ (with $F(v)=0$ for $v<1$). Then $\hat h (\omega , z)=0$ for  characters with non-trivial restriction to  $C_{\Q,1}=\hatz$, while  $\hat h (1 , z)=\Phi(z-\frac 12)$. Moreover note that
for the archimedean place $v$ of $\K=\Q$ one has, disregarding the principal values for simplicity, 
$$
\int_{\K_v^*} \frac{h(u^{-1}) }{   \vert 1-u   \vert} \, d^* u=\int_{\R^*} \frac{h(u) }{   \vert 1-u^{-1}   \vert} \, d^* u
$$
$$
=\frac 12\int_1^\infty h(t)\left( \frac{1}{\vert 1-t^{-1}\vert}+\frac{1}{\vert 1+t^{-1}\vert}\right)\frac{dt}{t}=\int_1^\infty\frac{t^{3/2}F(t)}{t(t^2-1)}dt
$$
where the $\frac 12$ comes from the normalization of the multiplicative Haar measure of $\R^*$ viewed as a modulated group.  In a similar way,  the normalization of the multiplicative Haar measure on $\Q_p^*$ shows that for the finite place associated to the prime $p$  one gets
 the term $\sum_{m=1}^\infty \log p \,\,p^{-m/2}F(p^m)$.

\subsection{RH for function fields}
\label{subsecff}

When the module $\mmod(\K)$ of a global field is a discrete subgroup of $ \R_+^*$ it is of the form $\mmod(\K)=q^\Z$ where $q$ is a prime power, and the field $\K$ is the function field of a smooth projective curve $C$ over the finite field $\F_q$. 

Already at the beginning of the XX-th century, Emil Artin and Friedrich Karl Schmidt have generalized RH to the case of function fields. 
We refer to the text of Cartier \cite{Cart}  where he explains how Weil' \!\!s definition of the zeta function associated to a variety over a finite field slowly emerged, starting with the thesis of E. Artin where this zeta function was defined for quadratic extensions of $\F_q[T]$, explaining F. K. Schmidt' \!\!s generalization to arbitrary extensions of $\F_q[T]$ and the work of Hasse on the ``Riemann hypothesis" for elliptic curves over finite fields.

 When the global field $\K$ is a function field, geometry comes to the rescue.  The problem becomes intimately related to the geometric one of estimating the number $N(q^r):=\#\,C(\F_{q^r})$ of points of $C$ rational over a finite extension $\F_{q^r}$ of the field of definition of $C$. The analogue of the Riemann zeta function is a generating function: the Hasse-Weil zeta function   
\begin{equation}\label{HW}
\zeta_C(s):=Z(C,q^{-s}), \ \  Z(C,T) := \exp\left(\sum_{r\geq 1}N(q^r)\frac{T^r}{r}\right)
\end{equation}

The analogue of RH for $\zeta_C$ was proved by  Andr\' e Weil in 1940. Pressed by the circumstances (he was detained in jail) he sent a Comptes-Rendus note to E. Cartan announcing his result. Friedrich Karl Schmidt  and Helmut Hasse had previously been able to transpose the Riemann-Roch theorem in the framework of geometry over finite fields and shown its implications for the zeta function: it is a rational fraction (of the variable $T$) and it satisfies a functional equation. But it took Andr\'e Weil several years  to put on solid ground a general theory of algebraic geometry in finite characteristic that would justify his geometric arguments and  allow him to 
transpose the Hodge index theorem in the form due to the Italian geometers  Francesco Severi and Guido Castelnuovo at the beginning of the  XX-th century.

\subsection{The proof using Riemann-Roch on $\bar C\times \bar C$}\label{rrstrat}

Let  $C$ be  a smooth projective curve  over the finite field $\F_q$. 
The first step is to extend the scalars from $\F_q$ to an algebraic closure $\bar\F_q$. Thus one lets 
\begin{equation}\label{extscal}
\bar C:=C\otimes_{\F_q}\bar\F_q
\end{equation}
This operation of extension of scalars does not change the points over $\bar\F_q$, \ie one has $\bar C(\bar\F_q)=C(\bar\F_q)$. The Galois action of the Frobenius automorphism of $\bar\F_q$ raises the coordinates of any point $x\in C(\bar\F_q)$ to the $q$-th power and this transformation of $C(\bar\F_q)$ coincides with the {\em relative Frobenius} $\fr_r:=\fr_C\times \id$ of $\bar C$, where $\fr_C$ is the {\em absolute Frobenius} of $C$ (which is the identity on points of the scheme and the $q$-th power map in the structure sheaf). The relative Frobenius $\fr_r$ is $\bar\F_q$-linear by construction and one can consider its graph 
in the   surface $X=\bar C\times_{\bar\F_q} \bar C$ which is the square of $\bar C$. This graph is the Frobenius correspondence $\Psi$.   It is important to work over an algebraically closed field in order to have a good intersection theory. 
This allows one to express the right hand side of the explicit formula \eqref{weil4} for the zeta function $\zeta_C$ as an intersection number $D.\Delta$, where $\Delta$ is the diagonal in the square and $D=\sum a_k\Psi^k$ is the divisor given by a  finite integral linear combination of  powers of the Frobenius correspondence. The terms $\hat h (0)$, $\hat h (1)$ in the explicit formula are also given by intersection numbers $D.\xi_j$, where  
\begin{equation}\label{xin}
   \xi_0=e_0\times \bar C\,, \ \xi_1=\bar C\times e_1
\end{equation}
where the $e_j$ are points of $\bar C$. One then  
 considers divisors on $X$  up to
the additive subgroup of principal divisors \ie those corresponding to an element  $f\in \cK$ of the function field of $X$.
 The problem is then reduced to proving the negativity of $D.D$ (the self-intersection pairing) for divisors of degree zero. The Riemann-Roch theorem on the surface $X$ gives the answer. 
 To each divisor $D$ on $X$ corresponds an index problem and one has a finite dimensional vector space of solutions $H^0(X,\cO(D))$ over $\bar\F_q$. Let
\begin{equation}\label{dim}
    \ell(D)={\rm dim}\, H^0(X,\cO(D))
\end{equation}
The best way to think of the sheaf $\cO(D)$ is in terms of Cartier divisors, \ie a global section of the quotient sheaf $\cK^\times/\cO_X^\times$, where $\cK$ is the constant sheaf corresponding to the function field  of $X$ and $\cO_X$ is the structure sheaf. The sheaf $\cO(D)$ associated to a Cartier divisor is obtained by taking the sub-sheaf of $\cK$ whose sections on $U_i$ form the sub $\cO_X$-module generated by $f_i^{-1}\in \Gamma(U_i,\cK^\times)$ where the $f_i$ represent $D$ locally. One has a ``canonical" divisor $K$ and Serre duality
\begin{equation}\label{Sdual}
    {\rm dim}\, H^2(X,\cO(D))={\rm dim}\, H^0(X,\cO(K-D))
\end{equation}
Moreover the following Riemann-Roch formula holds
\begin{equation}\label{RRform}
    \sum_0^2(-1)^j {\rm dim}\, H^j(X,\cO(D))=\frac 12 D.(D-K)+\chi(X)
\end{equation}
where $\chi(X)$ is the arithmetic genus.
All this yields the Riemann-Roch  inequality
\begin{equation}\label{rrine}
    \ell(D)+ \ell(K-D)\geq \frac 12 D.(D-K)+\chi(X)
\end{equation}
One then applies Lemma \ref{simplelem1} to the quadratic form $\inter(D,D')=D.D'$ using the
$\xi_j$ of \eqref{xin}. One needs three basic facts (\cite{grmt})
\begin{enumerate}
  \item If $\ell(D)>1$ then $D$ is equivalent to a strictly positive divisor.
  \item If $D$ is a strictly positive divisor then
  $$ D.\xi_0+D. \xi_1>0$$
  \item One has $\xi_0.\xi_1=1$ and $\xi_j.\xi_j=0$.
\end{enumerate}
One then uses \eqref{rrine} to show (see \cite{grmt}) that if $D.D>0$ then after a suitable rescaling by $n>0$ or $n<0$ one gets $\ell(nD)>1$ which shows that the hypothesis (2) of the following simple Lemma \ref{simplelem1} is fulfilled, and hence that RH holds for $\zeta_C$,
 \begin{lem} \label{simplelem1} Let $\inter(x,y)$ be a symmetric bilinear form on a  vector space $E$ (over $\Q$ or $\R$). Let $\xi_j \in E$, $j\in \{0,1\}$, be such that
\begin{enumerate}
  \item $\inter(\xi_j,\xi_j)=0$ and $\inter(\xi_0,\xi_1)=1$.
  \item For any $x\in E$ such that $\inter(x,x)>0$ one has $\inter(x,\xi_0)\neq 0$ or $\inter(x,\xi_1)\neq 0$.
 \end{enumerate}
 Then  one has the inequality
 \begin{equation}\label{negative4}
    \inter(x,x)\leq 2 \inter(x,\xi_0)\inter(x,\xi_1)\qqq x\in E
 \end{equation}
\end{lem}
The proof takes one line but the meaning of this lemma is to reconcile the ``naive positivity" of the right hand side of the explicit formula \eqref{weil4} (which is positive when $h\geq 0$ vanishes near $u=1$) with the negativity of the left hand side needed to prove RH (\cff\S \ref{sectSelberg} \eqref{negcrit} below).

\vsp
\centerline{\colorbox{lightgray}{\parbox[top][3cm][c]{13cm}{
At this point we see that it is highly desirable to find a geometric framework 
for the Riemann zeta function itself, in which the Hasse-Weil formula \eqref{HW}, the geometric interpretation of the explicit formulas, the Frobenius correspondences, the divisors, principal divisors, Riemann-Roch problem on the curve and the square of the curve all make sense. 
}}}

\vsp

Such a tentative framework will be explained in Section \ref{algeomattack}. It involves in particular the refinement of the notion of geometric space which was uncovered by Grothendieck and to which we now briefly turn.

\subsection{Grothendieck and the notion of topos}
\label{subsecff}

The essential ingredients of the proof explained in \S \ref{rrstrat} are the intersection theory for divisors on $\bar C\times \bar C$, sheaf cohomology and Serre duality, which give the formulation of the Riemann-Roch theorem. Both  owe to the discovery of sheaf theory by J. Leray and the pioneering work of J. P. Serre on the use of sheaves for the Zariski topology in the algebraic context, with his fundamental theorem comparing the algebraic and analytic frameworks. The next revolution came from the elaboration by A. Grothendieck and M. Artin of etale $\ell$-adic cohomology. It allows one to express the Weil zeta function of a smooth projective variety $X$ defined over a finite field $\F_q$ \ie the function $Z(X,t)$ given by \eqref{HW} with $t=q^{-s}$ which continues to make sense in general, as an alternate product of the form
\begin{equation}\label{etalecohomol}
Z(X,t)=\prod_{j=0}^{2\, \dim X} \det(1-t F^*\mid H^j({\bar X}_\ett,\Q_\ell))^{(-1)^{j+1}}
\end{equation}
where $F^*$ corresponds to the action of the Frobenius on the $\ell$-adic cohomology and $\ell$ is a prime which is prime to $q$. This equality follows from a Lefschetz  formula for the number $N(q^r)$ of fixed points of the $r$-th power of the Frobenius and when $X=C$ is a curve the explicit formulas reduce to the Lefschetz formula. The construction of the cohomology groups $H^j(\bar X_{\ett},\Q_\ell)$ is indirect and they are defined as :
$$
H^j(\bar X_{\ett},\Q_\ell)=\varprojlim_n \left( H^j(\bar X_\ett,\Z/\ell^n\Z) \right)\otimes_{\Z_\ell}\Q_\ell
$$
where $\bar X_{\ett}$ is the etale site of $\bar X$. Recently the etale site of a scheme has been refined \cite{BS} to the {\em pro-etale} site whose objects   no longer satisfy any  finiteness condition. The cohomology groups $H^j(\bar X_{\rm proet},\bar\Q_\ell)$ are then directly obtained using the naive interpretation  (without torsion coefficients). 
One needs to pay attention in \eqref{etalecohomol} to the precise definition of $F$, it is either the relative Frobenius $\fr_r$ or the {\em Geometric Frobenius} $\fr_g$ which is the inverse of the {\em Arithmetic Frobenius} $\fr_a$. The product $\fr_a\circ \fr_r=\fr_r\circ \fr_a$ is the absolute Frobenius $\fr$ which acts trivially on the $\ell$-adic cohomology. To understand the four different incarnations of ``the Frobenius" it is best to make them explicit in the simplest example of the scheme $\Spec R$ where $R=\bar\F_q[T]$ is the ring of polynomials $P(T)=\sum a_j T^j$, $a_j\in \bar\F_q$
\begin{itemize}
\item Geometric Frobenius: $\sum a_j T^j\mapsto \sum a_j^{1/q} T^j$
\item Relative Frobenius: $P(T)\mapsto P(T^q)$
\item  Absolute Frobenius: $P(T)\mapsto P(T)^q$
\item  Arithmetic Frobenius: $\sum a_j T^j\mapsto \sum a_j^q T^j$
\end{itemize}
The motivation of Grothendieck for developing etale cohomology came from the search of a Weil cohomology and the Weil conjectures which were solved by Deligne in 1973 (\cite{Deligne}).

 In his quest Grothendieck uncovered several key concepts such as those of schemes and above all that of topos, in his own words:

\begin{quote}{\em
C'est le th\`eme du
topos, et non celui des sch\'emas, qui est ce ``lit'', ou cette ``rivi\`ere profonde'', o\`u viennent
s'\'epouser la g\'eom\'etrie et l'alg\`ebre, la topologie et l'arithm\'etique, la logique math\'ematique
et la th\'eorie des cat\'egories, le monde du continu et celui des structures ``discontinues'' ou
``discr\`etes''. Si le th\`eme des sch\'emas est comme le {\it c\oe ur} de la g\'eom\'etrie nouvelle, le th\`eme
du topos en est l'enveloppe, ou la {\it demeure}. Il est ce que j'ai con\c cu de plus vaste, pour
saisir avec finesse, par un m\^eme langage riche en r\'esonances g\'eom\'etriques, une ``essence''
commune \`a des situations des plus \'eloign\'ees les unes des autres, provenant de telle r\'egion
ou de telle autre du vaste univers des choses math\'ematiques.
}\end{quote}

\section{Riemannian Geometry, Spectra  and trace formulas} \label{analysisattack}
Riemannian Geometry gives a wealth of ``spectra" of fundamental operators associated to a geometric space, such as the Laplacian and the Dirac operators.
\subsection{The Selberg trace formula}\label{sectSelberg}
In the case of compact Riemann surfaces $X$ with constant negative curvature $-1$, the Selberg trace formula \cite{Selberg}, takes the following form where the eigenvalues of the Laplacian are written in the form\footnote{where the argument of $r_n$ is either $0$ or $-\pi/2$} $\lambda_n=-(\frac 14+r_n^2)$. 
Let $\delta>0$, $h(r)$ be an analytic function in the strip $\vert \Im(r)\vert\leq \frac 12+\delta$ and such that $h(r)=h(-r)$ and with $(1+r^2)^{1+\delta}\vert h(r)\vert$ being bounded. Then \cite{Selberg, [Se],Hejhal}, with $A$ the area of $X$,
\begin{equation}\label{Selberg}
\sum h(r_n)=\frac{A}{4\pi}\int_{-\infty}^\infty {\rm tanh}(\pi r)h(r)rdr 
+\sum_{\{T\}}\frac{\log N(T_0)}{N(T)^\frac 12-N(T)^{-\frac 12}}g(\log N(T))
\end{equation}
where $g$ is the Fourier transform of $h$, \ie more precisely 
$g(s)=\frac{1}{2\pi}\int_{-\infty}^\infty h(r)e^{-irs}dr$. The $\log N(T)$ are the lengths of the periodic orbits of the geodesic flow with $\log N(T_0)$ being the length of the primitive one. 
Already in 1950-51, Selberg saw the striking similarity of his  formula with \eqref{explfor} which (\cf~\cite{Hejhal}) can be rewritten in the following form, with $h$ and $g$ as above and the non-trivial zeros of zeta expressed in the form 
 $\rho=\frac 12+i\gamma$,  
\begin{equation}\label{Hejhal}
\sum_\gamma h(\gamma)=h(\frac i2)+h(-\frac i2)+\frac{1}{2\pi}\int_{-\infty}^\infty \omega(r)h(r)dr-2\sum \Lambda(n)n^{-\frac 12}g(\log n)
\end{equation}
where
$$
\omega(r)=\frac{\Gamma'}{\Gamma}\left(\frac 14+i\frac r2\right)-\log \pi, \  \  \frac{\Gamma'}{\Gamma}\left(s\right)=\int_0^1\frac{1-t^{s-1}}{1-t}dt-\gamma\qqq s,\Re(s)>0
$$
and $\Lambda(n)$ is the von-Mangoldt function with value $\log p$ for powers $p^\ell$ of  primes and zero otherwise. 
Moreover Selberg found that there is a zeta function which corresponds to \eqref{Selberg} in the same way that $\zeta(s)$ corresponds to \eqref{Hejhal}.
The role of Hilbert space is crucial in the work of Selberg to ensure that the zeros of his zeta  function satisfy the analogue of RH. This role of Hilbert space is implicit as well in RH which has been reformulated by Weil as the positivity of the functional $W(g)$ defined as both sides  of \eqref{Hejhal}. More precisely the equivalent formulation  is that  $W(g\star g^*)\geq 0$ on functions $g$ which correspond to Fourier transforms of analytic functions $h$ as above (\ie even and analytic in a strip $\vert\Im z\vert \leq \frac 12+\delta$) where for even functions one has $g^*(s):=\overline{ g(-s)}=\overline{ g(s)}$. Moreover by \cite{B3,burnol1},  it is enough, using Li's criterion (\cff\cite{Li,B3}), to check the positivity on a small class of explicit real valued functions with compact support. In fact for later purposes it is better to write this criterion as
\begin{equation}\label{negcrit}
RH \iff \inter(f,f)\leq 0 \qqq f \mid \int f(u)d^*u=\int f(u)du=0
\end{equation}
where for real compactly supported functions on $\R_+^*$, we let $\inter(f,g):=N(f\star \tilde g)$ where $\star$ is the convolution product on $\R_+^*$,  $\tilde g(u):=u^{-1}g(u^{-1})$, and
\begin{equation}\label{negcrit1}
N(h):= \sum_{n=1}^\infty \Lambda(n)h(n)+ \int_1^\infty\frac{u^2h(u)-h(1)}{u^2-1}d^*u+c\, h(1)\,, \ c=\frac12(\log\pi+\gamma)
\end{equation}

The Selberg trace formula has been considerably extended by J. Arthur and plays a key role in the Langland's program. We refer to  \cite{Arthur} for an introduction to this vast topic.

\subsection{The minus sign and absorption spectra}
\begin{figure}[H]
\hfill\begin{minipage}[c]{0.38\textwidth}
    \caption{ \small{\it The three kinds of spectra occuring in spectroscopy: 
1) The top one is the ``continuous spectrum" which occurs when white light is decomposed by passing through a prism. 2) The middle one is the ``emission spectrum" which occurs when the light emitted by a heated gas is decomposed by passing through a prism and gives shining lines-a signature of the gas-over a dark background. 3) The third one is the ``absorption spectrum" which occurs when white light traverses a cold gas and is then decomposed by passing through a prism. It appears as dark lines in a background continuous spectrum. The absorption lines occur at the same place as the emission lines.
       }} \label{absspectrum}
  \end{minipage}\hfill
\begin{minipage}[c]{0.08\textwidth}  \end{minipage}
  \begin{minipage}[c]{0.48\textwidth}
    \includegraphics[width=\textwidth]{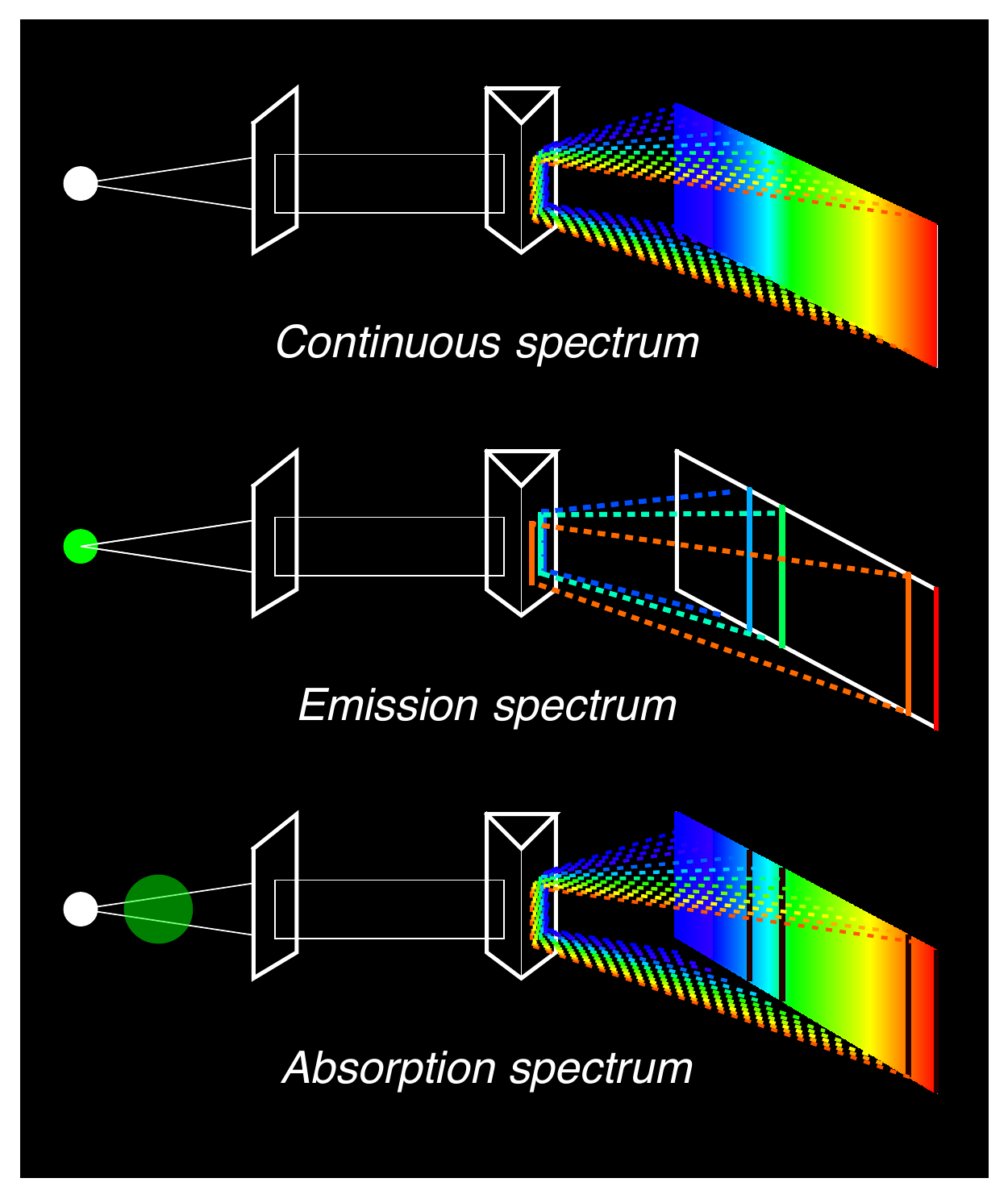}
  \end{minipage}  
\end{figure}

The Selberg  trace formula \cite{Selberg, [Se]} for Riemann surfaces of finite area, acquires additional terms which make it look \eg in the case of $X=H/PSL(2,\Z)$ (where $H$ is the upper half plane with the Poincar\'e metric) even more similar to the explicit formulas, since the parabolic terms now involve explicitly the sum $$2\sum_{n=1}^\infty \frac{\Lambda(n)}{n}g(2\log n)$$  
Besides the square root in  the $\Lambda(n)$ terms in the explicit formulas \eqref{Hejhal}
$$
-2\sum_{n=1}^\infty \frac{\Lambda(n)}{n^{\frac 12}}g(\log n)
$$
there is however a striking difference which is that these terms occur with a positive sign instead  of the   negative sign in  \eqref{Hejhal}, as discussed in \cite{Hejhal} \S 12. This discussion of the minus sign was extended to the case of the semiclassical limit of Hamiltonian
systems in physics in \cite{Berry}. In order to get some intuition of what this reveals, it is relevant to go back to the origin of spectra in physics, \ie to the very beginning of spectroscopy. It occurred  when Joseph Von Fraunhofer (1787-1826)  could identify, using self-designed instruments, about 500 dark lines in the light coming from the sun, decomposed using the dispersive power of a spectroscope such as a prism (\cff Figure \ref{absspectrum}). These dark lines constitute the ``absorption spectrum" and it took about 45 years  before  Kirchhoff and Bunsen noticed that several of these Fraunhofer lines coincide (\ie have the same wave length) with the bright lines of the  ``emission" spectrum of heated elements, and showed that they could be reobtained by letting white light traverse a cold gas.
In his work on the trace formula in the finite covolume case, Selberg had to take care of a superposed continuous spectrum due to the presence of the non-compact cusps of the Riemann surface.

\subsection{The adele class space and the explicit formulas}\label{adeleclass}

I had the chance to be invited at the Seattle meeting in 1996 for the celebration of the proof of the prime number theorem. The reason was the paper \cite{BC} (inspired from \cite{[J2]}) in which the Riemann zeta function appeared naturally as the partition function of a quantum mechanical system (BC system) exhibiting phase transitions. The RH had been at the center of discussions in the meeting and I knew the analogy between the BC-system and the set-up that V. Guillemin proposed in \cite{gui} to explain the Selberg trace formula using the action of the geodesic flow on the horocycle foliation. To a foliation is associated a von Neumann algebra \cite{Co-foliations}, and the horocycle foliation on the sphere bundle of a compact Riemann surface gives a factor of type II$_\infty$ on which the geodesic flow acts by scaling the trace. An entirely similar situation comes canonically from the BC-system at critical temperature and after interpreting the dual system in terms of adeles, I was led by this analogy to consider the action of the idele class group of $\Q$ on the adele class space, \ie the quotient $\Q^\times\backslash \A_\Q$ of the adeles $\A_\Q$ of $\Q$ by the action of $\Q^\times$. I knew from the BC-system that the action of $\Q^\times$, which preserves the additive Haar measure, is ergodic for this measure and gives the same factor of type II$_\infty$ as the horocycle foliation. Moreover the dual action scales the trace in the same manner. 

Let $\K$ be a global field and $C_\K=\GL_1(\A_\K)/\K^\times$ the idele class group. The module $\mmod:C_\K\to \R_+^*$ being proper with cocompact range, one sees that the Haar measure on the Pontrjagin dual group of $C_\K$ is diffuse. Since a point is of measure $0$ in a diffuse measure space there is no way one can see the absorption spectrum without introducing some smoothness on this dual which is done using a Sobolev space $L_{\delta  }^2 ( C_\K )$ of functions on $C_\K$ which (for fixed $\delta>1$) is defined as
\begin{equation}
||\xi||^2=\int_{C_\K} \vert \xi (x) \vert^2 \, \rho(x) \, d^* x, \  \,  \ \ \rho(x):=(1 + \log \vert x
\vert^2)^{\delta   / 2} \label{sobolevck}
\end{equation}
\begin{defn}\label{adclassspace}
 Let $\K$ be a global field, the adele class space of $\K$ is the quotient $X_\K=\A_\K/\K^\times$ of the adeles of $\K$ by the action of $\K^\times$ by multiplication. 
\end{defn}
We then  
consider the codimension 2 subspace $\cS
(\A_\K)_0$ of the Bruhat-Schwartz space
$\cS (\A_\K)$ (\cf~\cite{Bruhat}) given by the conditions
$
f(0) = 0 \, , \ \int f \, dx = 0 
$
The Sobolev space $L_{\delta  }^2 (X_\K)_0$ is the separated completion of $\cS
(\A_\K)_0$ for the norm with square 
\begin{equation}
||f||^2=\int_{C_\K} \vert \sum_{q\in \K^*}  f (qx) \vert^2 \, \rho(x) \, \vert x\vert d^* x \label{spec2}
\end{equation}
Note that by construction all functions of the form $f(x)=g(x)-g(qx)$ for some $q\in \K^\times$ belong to the radical of the norm \eqref{spec2}, which corresponds to the operation of quotient of Definition \ref{adclassspace}. In particular  the representation of ideles on $\cS (\A_\K)$
given by 
\begin{equation}
(\urep (\alpha) \xi) (x) = \xi (\alpha^{-1} x)  \  \   \forall \,
\alpha \in \GL_1(\A_\K) \, , \ x \in \A_\K \,  \label{spec(3)}
\end{equation}
induces a representation $\urep_a$ of $C_\K$ on $L_{\delta  }^2 (X_\K)_0$. 
One has by construction a natural isometry $\mapE:L_{\delta  }^2 (X_\K)_0\to L_{\delta  }^2 ( C_\K )$ which intertwines the representation $\urep_a$ with the regular representation of $C_\K$ in $L_{\delta  }^2 ( C_\K )$ multiplied by the square root of the module.  This representation restricts to the cokernel of the map $\mapE$, which splits as a direct sum of subspaces labeled by the characters  of the compact group $C_{\K,1}=\Ker\, \mmod$ and its spectrum in each sector gives the  zeros of $L$-functions with Gr\"ossencharakter. The shortcoming of this construction is in the artificial weight $\rho(x)$, which is needed to see this absorption spectrum  but  only sees the zeros which are on the critical line and where the value of $\delta$ artificially cuts the multiplicities of the zeros (\cf~\cite{Co-zeta}). 
 
This state of affairs is greatly improved if one gives up trying to prove RH but retreats to an interpretation of the explicit formulas as a trace formula. One simply replaces the above Hilbert space set-up by a softer one involving nuclear spaces \cite{Meyer}. The spectral side now involves all non-trivial zeros and, using the preliminary results of \cite{Co-zeta,Co99,burnol} one gets that the geometric side is given by:
\begin{equation}\label{geomside}
\Tr_{\rm distr}\left(\int h(w)\urep(w)d^*w\right )=\sum_v\int_{\K^\times_v}\,\frac{h(w^{-1})}{|1-w|}\,d^*w
\end{equation}
We refer to \cite{Co-zeta, Meyer, CMbook} for a detailed treatment.
The subgroups $\K^\times_v\subset C_\K=\GL_1(\A_\K)/\GL_1(\K)$ arise as isotropy groups.
One can understand why the terms $\displaystyle \frac{h(w^{-1})}{|1-w|}$ occur in the trace formula by computing, formally as follows, the trace  of the scaling operator $T=\urep_{w^{-1}}$ when working on the local field $\K_v$ completion of the global field $\K$ at the place $v$, one has
$$T\xi(x)=\xi(w x)=\int k(x,y)\xi(y)dy\,
 $$
so that $T$ is given by the distribution kernel $k(x,y)=\delta(w x-y)$ and its trace is 
 $$
\Tr_{\rm distr}(T)=\int k(x,x)\,dx=\int \delta(w x-x)\,dx=\frac{1}{|w-1|}\int \delta(z)\,dz=\frac{1}{|w-1|}
$$

When working at the level of adeles one treats all places on the same footing and thus there is an overall minus sign in front of the spectral contribution. 
Thus the Riemann spectrum appears naturally as an absorption spectrum from the adele class space. As such, it is difficult to show that it is ``real". While this solves the problem of giving a trace formula interpretation of the explicit formulas, there is of course still room for an interpretation as an emission spectrum. However from the adelic point of view it is unnatural to separate the contribution of the archimedean place.

\section{The Riemann-Roch strategy: A Geometric Framework}
\label{algeomattack}
In this section we shall present a geometric framework which has emerged over the years in our joint work with C. Consani and seems suitable in order to transpose the geometric proof of Weil to the case of RH. The aim is to apply the Riemann-Roch strategy of \S \ref{rrstrat}. The geometry involved will be of elaborate  nature inasmuch as it relies on the following three theories:
\begin{enumerate}
\item{Noncommutative Geometry.}
\item{Grothendieck topoi.}
\item{Tropical Geometry.}
\end{enumerate}

\subsection{The limit $q\to 1$ and the Hasse-Weil formula} \label{counting}

In \cite{Soule} (\cf~\S 6), C. Soul\'e, motivated by \cite{Man-zetas} (\cf~\S 1.5) and  \cite{Steinberg, Tits, Ku, Den1,Den2, Kapranov}, introduced the zeta function of a variety $X$ over $\F_1$ using the {\em polynomial} counting function $N(x)\in\Z[x]$ associated to $X$. The definition of the zeta function is as follows
\begin{equation}\label{zetadefn}
\zeta_X(s):=\lim_{q\to 1}Z(X,q^{-s}) (q-1)^{N(1)},\qquad s\in\R
\end{equation}
where $Z(X,q^{-s})$ denotes the evaluation at $T=q^{-s}$  of the Hasse-Weil exponential series
\begin{equation}\label{zetadefn1}
Z(X,T) := \exp\left(\sum_{r\ge 1}N(q^r)\frac{T^r}{r}\right)
\end{equation}
For instance, for a projective space $\P^n$ one has $N(q)=1+q+\ldots +q^n$ and 
$$
 \zeta_{\P^n(\F_1)}(s) =\lim_{q\to 1} (q-1)^{n+1} \zeta_{\P^n(\F_q)}(s) = \frac{1}{\prod_0^n(s-k)}
$$
It is natural to wonder on the existence of a ``curve'' $C$ suitably defined over $\F_1$, whose zeta function $\zeta_C(s)$ is the complete Riemann zeta function  $\zeta_\Q(s)=\pi^{-s/2}\Gamma(s/2)\zeta(s)$ (\cf~also \cite{Man-zetas}). The first step is to find a counting function $N(q)$ defined for $q\in [1,\infty)$ and such that  \eqref{zetadefn} gives $\zeta_\Q(s)$. But there is an obvious difficulty since as $N(1)$ represents the Euler characteristic one should expect that   $N(1)=-\infty$ (since the dimension of $H^1$ is infinite). This precludes the use of \eqref{zetadefn}  and also seems to contradict the expectation that $N(q)\geq 0$ for $q\in (1,\infty)$. As shown in \cite{CC0,CC1} there is a simple way to solve the first difficulty by passing to the logarithmic derivatives of both terms in equation \eqref{zetadefn} and observing that the Riemann sums of an integral appear from the right hand side. One then gets instead of  \eqref{zetadefn} the equation:
\begin{equation}\label{logzetabis}
    \frac{\partial_s\zeta_N(s)}{\zeta_N(s)}=-\int_1^\infty  N(u)\, u^{-s}d^*u
\end{equation} 
Thus the integral equation \eqref{logzetabis} produces a precise equation for the counting function $N_C(q)=N(q)$ associated to $C$:
\begin{equation}\label{special}
   \frac{\partial_s\zeta_\Q(s)}{\zeta_\Q(s)}=-\int_1^\infty  N(u)\, u^{-s}d^*u
\end{equation}
One finds that this equation admits a solution which is a {\em distribution} and is given with  $ \varphi(u):=\sum_{n<u}n\,\Lambda(n)$, by the equality
\begin{equation}\label{Nu}
    N(u)=\frac{d}{du}\varphi(u)+ \kappa(u)
\end{equation}
where  $\kappa(u)$ is the distribution
which appears in the explicit formula \eqref{explfor},
$$
\int_1^\infty\kappa(u)f(u)d^*u=\int_1^\infty\frac{u^2f(u)-f(1)}{u^2-1}d^*u+cf(1)\,, \qquad c=\frac12(\log\pi+\gamma)
$$
The conclusion is that the distribution $N(u)$ is positive on $(1,\infty)$ and  is given  by
\begin{equation}\label{fin2}
    N(u)=u-\frac{d}{du}\left(\sum_{\rho\in Z}{\rm order}(\rho)\frac{u^{\rho+1}}{\rho+1}\right)+1
\end{equation}
where the derivative is taken in the sense of distributions, and the value at $u=1$ of the  term
 $\displaystyle{\omega(u)=\sum_{\rho\in Z}{\rm order}(\rho)\frac{u^{\rho+1}}{\rho+1}}$ is given  by
$\frac 12+ \frac \gamma 2+\frac{\log4\pi}{2}-\frac{\zeta'(-1)}{\zeta(-1)}
$.

\begin{figure}[H]
\hfill\begin{minipage}[c]{0.27\textwidth}
    \caption{ \small{\it This represents a function $J(u)$ which is a primitive of the counting distribution $N(u)$. This function is increasing and tends to $-\infty$ when $u\to 1$. The wiggly graph represents the approximation of $J(u)$ obtained using the symmetric set $Z_m$ of the first $2m$ zeros, by $J_m(u)=\frac{u^2}{2}-\sum_{Z_m}{\rm order}(\rho)\frac{u^{\rho+1}}{\rho+1}+u$ Note that $J(u)\to -\infty$ when $u\to 1+$.
       }} \label{figcounting} 
  \end{minipage}\hfill
  \begin{minipage}[c]{0.60\textwidth}
    \includegraphics[width=\textwidth]{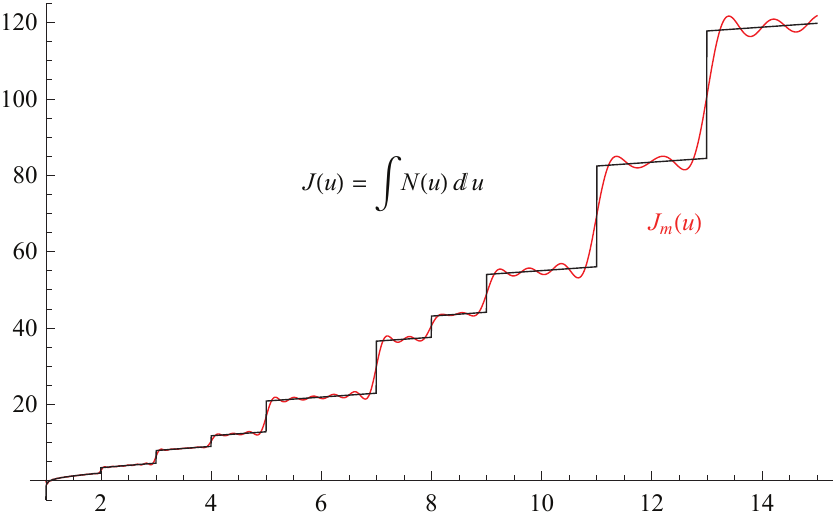}
  \end{minipage}  
\end{figure}

The primitive $J(u)=\frac{u^2}{2}-\omega(u)+u$ of $N(u)$ is an increasing function on $(1,\infty)$, but tends to $-\infty$ when $u\to 1+$ while its value $J(1)$ is finite. The tension between the positivity of the distribution $N(q)$ for $q>1$ and the expectation that its  value $N(1)$ should be $N(1)=-\infty$ is resolved by the theory of distributions: $N$ is {\em finite} as a distribution, but when one looks at it as a function its value at $q=1$ is formally given by
$$
N(1)=2-\lim_{\epsilon\to 0}\frac{\omega(1+\epsilon)-\omega(1)}{\epsilon}\sim-\frac 12 E \log E,\qquad \ E=\frac 1\epsilon
$$
which is $-\infty$ and in fact reflects, when $\epsilon\to 0$,  the density of the zeros. Note that this holds independently of the choice of the principal value in the explicit formulas.  This subtlety does not occur for function fields $\K$ since their module $\mmod(\K)$ is discrete so that distributions and functions are the same thing. There is one more crucial nuance between the case $\K=\Q$ and the function fields: the distribution $\kappa(u)$ which is the archimedean contribution to $N(u)$ in \eqref{Nu}, does not fulfill the natural inequality $N(q)\leq N(q^r)$ expected of a counting function. This is due to the  terms $\vert 1-u\vert^{-1}$ in the Weil explicit formula, which as explained in \S \ref{sectweilexpl} contribute non-trivially at the archimedean place, and indicate that the counting needs to take into account an ambient larger space and transversality factors as in \cite{gui}.
In fact, we have seen in Section \ref{adeleclass} that the noncommutative space of adele classes of a global field provides a framework to interpret the explicit formulas of Riemann-Weil in number theory as a trace formula, and that the geometric contributions give the right answer. In \cite{CC1}, we showed that the quotient
\begin{equation}\label{doublequot}
X_\Q:=\Q^\times\backslash \A_\Q/\hatz
\end{equation}
of the adele class space $\Q^\times\backslash \A_\Q$ of the rational numbers by the maximal compact subgroup $\hatz$ of the idele class group,  gives by considering the induced action of $\R_+^\times$, the above counting distribution $N(u)$, $u\in [1,\infty)$, which determines, using the Hasse-Weil formula in the limit $q\to 1$, the complete Riemann zeta function. The next step is to understand that  the action of $\R_+^\times$ on the space $X_\Q$ is in fact the action of the Frobenius automorphisms $\fr_\lambda$ on the points of the arithmetic site-- an object of algebraic geometry--over $\rmax$. To explain this we first need to take an excursion in the exotic world of ``characteristic one".

\subsection{The world of characteristic $1$}\label{sectchar1}
The key words here are: Newton polygons, Thermodynamics, Legendre transform, Game theory, Optimization, Dequantization, Tropical geometry. One alters the basic operation of addition of positive real numbers, replacing $x+y$ by  $x\vee y:=\max(x,y)$. When endowed with this operation as addition and with the usual multiplication, the positive real numbers become a semifield $\rmax$. It is of characteristic $1$, \ie $1\vee 1=1$ and contains the smallest semifield of characteristic $1$, namely the Boolean semifield $\B=\{0,1\}$. Moreover, $\rmax$ admits non-trivial automorphisms and one has 
$$
{\rm Gal}_\B(\rmax):=\Aut_\B(\rmax)=\R_+^*, \ \ \fr_\lambda(x)=x^\lambda \qqq x\in \rmax, \ \lambda \in \R_+^*
$$
thus providing a first glimpse of an answer to Weil's query in \cite{Weilcdc} of an algebraic framework in which the connected component of the idele class group would appear as a Galois group. More generally, for any abelian ordered group $H$ we let $H_{\rm max}=H\cup \{-\infty\}$ be the semifield obtained from $H$ by the max-plus construction, \ie the addition is given by the max, and the multiplication by $+$. In particular $\R_{\rm max}$ is isomorphic to $\rmax$ by the exponential map (\cf~\cite{Gaubert}). Historically, and besides the uses of $\rma$ in idempotent analysis and tropical geometry which are discussed below, an early use of $\rma$ occurred in the late fifties in the work of R. Cuninghame-Green in Birmingham, who established the spectral theory of irreducible matrices with entries in $\rma$ (\cff\cite{Cunni}) and in the sixties, in Leningrad, where Vorobyev used the $\rma$ formalism in his work motivated by combinatorial optimization, and proved a fundamental covering theorem. A systematic use of the $\rma$ algebra was developed by the  INRIA group at the beginning of the 80's  in their work on the modelization of discrete event systems \cite{Gaubert2}. We refer to \cite{Gaubert, Gaubert1} for a more detailed history of the subject, and for overwhelming evidence of its relevance in mathematics. We shall just give here a sample of this evidence starting by a really early occurrence in the work of C.G.J.~Jacobi\footnote{I am grateful to S. Gaubert for pointing out this early occurrence} and hoping to convince the reader that it would be a mistake  to dismiss this algebraic formalism and the analogy with ordinary algebra as trivial. 

\subsubsection{Optimization, Jacobi}
One of the early instances, around 1840, of the use of matrices over $\R_{\rm max}$ is the work of C.G.J.~Jacobi \cite{Jacobi} on optimal assignment problems, where he states 
\begin{quote}{\it \centerline{Problema}
Disponantur nn quantitates $h_k^{(i)}$ quaecunque in schema Quadrati, ita ut habeantur n series horizontales et n series verticales, quarum quaeque est n terminorum. Ex illis quantitatibus eligantur n transversales, \ie in seriebus horizontalibus simul atque verticalibus diversis positae, quod fieri potest n! modis; ex omnibus illis modis quaerendus est is, qui summam n numerorum electorum suppeditet maximam.
}\end{quote}
In other words, given a square matrix $m_{ik}=h_k^{(i)}$ he looks for the maximum over all  permutations $\sigma$  
of the quantity $\sum m_{j\sigma(j)}$. Using the algebraic rules of $\rma$ one checks that he is in fact computing the analogue of the determinant for the matrix $m_{ik}$. In fact the perfect definition of the determinant is more subtle and was obtained   
 in the work of Gondran-Minoux \cite{GM1}, instead of $\max\sum m_{j\sigma(j)}$ where $\sigma$ runs over all permutations,  one uses the signature of permutations and considers the pair $$({\rm det}_+(m_{ik}),{\rm det}_-(m_{ik})), \  \  {\rm det}_\pm(m_{ik})=\max \sum_{{\rm sign}(\sigma)=\pm} m_{j\sigma(j)}$$  The remarkable fact is that the Cayley-Hamilton theorem now holds, as the equality of two terms $P_+(m)=P_-(m)$ corresponding to the characteristic polynomial $P=(P_+,P_-)$. Each of the terms $P_\pm(m)\in M_n(\rma)$ is computed from the original matrix $m\in M_n(\rma)$ using the rules of matrices with entries in $\rma$ which turn $M_n(\rma)$ into a semiring.

\subsubsection{Idempotent analysis}

The essence of the theory of semiclassical analysis in physics rests in the comparison of quantum systems with their semiclassical counterpart, \cite{GS,Gutz,Gutz1,Maslov1,Berry}. In the eighties  V. P. Maslov and his collaborators developed a satisfactory algebraic framework which encodes the semiclassical limit of quantum mechanics. They called it idempotent analysis. We refer to  \cite{Maslov,Lit} for a detailed account and just mention briefly some salient features here.  The source of the variational formulations of  mechanics in the classical limit  is the behavior of  sums of exponentials
 $$
\sum e^{-\frac{ S_j}{ \hbar}} \sim e^{-\frac{ \inf S_j}{ \hbar}}, \  \ \text{when}\ \ \hbar\to 0
$$ 
which are, when $\hbar\to 0$, dominated by the contribution of the minimum of $S$. The starting observation is that one can encode this fundamental principle by simply conjugating the addition of numbers by the power operation $x\mapsto x^\epsilon$ and passing to the limit when $\epsilon\to 0$.   The new addition of positive real numbers is 
$$
\lim_{\epsilon \to 0}\left(x^{\frac 1\epsilon}+y^{\frac 1\epsilon}\right)^\epsilon=\max \{x,y\}=x\vee y
$$
and one recovers $\rmax$ as the natural home for semiclassical analysis.
The  superposition principle of quantum mechanics, \ie addition of vectors in Hilbert space,  now makes sense  in the limit and moreover the ``fixed point argument" proof of the Perron-Frobenius theorem works over $\rmax$ and shows that irreducible compact operators have one and only one eigenvalue\footnote{as mentioned above, this result was obtained already for matrices in 1962 by R. Cuninghame-Green}, thus reconciling classical determinism with the quantum variability. But the most striking discovery of this school of Maslov, Kolokolstov and Litvinov \cite{Maslov,Lit} is that the Legendre transform which plays a fundamental role in all of physics and in particular in 
 thermodynamics in the nineteenth century, is simply the Fourier transform  in the framework of  idempotent analysis!

  The contact between the INRIA school and the Maslov school was established in 92 when Maslov was invited  in the Seminar of Jacques Louis Lions in College de France. At the BRIMS HP-Labs workshop on Idempotency in Bristol (1994) organized by J. Gunawardena, several of the early groups of researchers in the field were there, and an animated discussion took place on how the field should be named.
The names ``max-plus", ``exotic", ``tropical", ``idempotent" were considered, each one having its defaults.

\subsubsection{Tropical geometry, Riemann-Roch theorems and the chip firing game}

The tropical semiring $\N_{\rm min}=\N\cup\{\infty\}$ with the operations $\min$ and $+$ was introduced by Imre Simon in \cite{Simon} to solve a decidability problem in rational language theory. His work is at the origin of the term ``tropical" used in tropical geometry which is a vast subject, see \eg \cite{Gelfand, Kap1, Mik,Sturm}. We refer to \cite{virotagaki} for an excellent introduction starting from the sixteenth Hilbert problem.  In its simplest form (\cff\cite{GK}) a tropical curve is given by a metric graph $\Gamma$ (\ie a graph with a usual line metric on its edges). The natural structure sheaf on $\Gamma$ is the sheaf $\cO$ of real valued functions which are continuous, convex, piecewise affine with integral slopes. The operations on such functions are given by the pointwise operations of $\rma$-valued functions, \ie $(f\vee g)(x)=f(x)\vee g(x)$ for all $x\in \Gamma$ and similar for the product which is given by pointwise addition. One also adjoins the constant $-\infty$ which plays the role of the zero element in the semirings of sections. One proceeds as in the classical case with the construction of the sheaf $\cK$ of semifields of quotients and finds the same type of functions as above but no longer convex. Cartier divisors make sense and one finds that the order of a section $f$ of $\cK$ at a point $x\in\Gamma$ is given by the sum of the (integer valued) outgoing slopes. The conceptual explanation of why the discontinuities of the derivative should be interpreted as zeros or poles is due to Viro, \cite{viro} who showed that it follows automatically if one understands that\footnote{as seen when using $\rma$ as the target of a valuation} the sum $x\vee x$ of two equal terms in $\rma$ should be viewed as ambiguous with all values in the interval $[-\infty, x]$ on equal footing. In their work Baker and Norine \cite{BN} proved in the discrete set-up of graphs (where $g$ is the genus and $K$ the canonical divisor) the
Riemann-Roch  equality in the form 
\begin{equation}\label{rr1}
r(D)-r(K-D)=\Deg(D)-g+1
\end{equation}
where by definition 
$r(D):=\max \{k\mid H^0(D-\tau )\neq \{-\infty\}\qqq \tau \geq 0,\ \Deg(\tau)=k\}$
and $H^0(D)$ is the $\rma$-module of global sections $f$ of the associated sheaf $\cO_D$ \ie sections of $\cK$ such that $D+(f)\geq 0$. 
The essence of the proof of \cite{BN} is that  the inequality ${\rm Deg}(D)\geq g$ for a divisor  implies $H^0(D)\neq \{-\infty\}$. Once translated in the language of the chip firing game (\opcit\!\!),  this fact is equivalent  to the existence of a winning strategy if one assumes that the total sum of dollars attributed to the vertices of the graph is $\geq g$ where $g$ is the genus. We refer to \cite{GK,MZ} for variants of the above Riemann-Roch  theorem, and to \cite{Dhar, Shor, PostS} for early occurrences of these ideas in a different context (including sandpile models and parking functions!).

\subsection{The arithmetic and scaling sites}

\subsubsection{The arithmetic site and Frobenius correspondences}

The  {\em arithmetic site} \cite{CC,CCas1} is an object  of algebraic geometry involving two elaborate mathematical concepts: the notion of topos and of (structures of) characteristic $1$ in algebra. 
A nice fact (\cf~\cite{Golan}) in characteristic $1$ is that, provided the semiring $R$ is \mc (\ie equivalently if it injects in its semifield of fractions) the map $x\mapsto x^n=\fr_n(x)$ is, for any integer $n\in \nt$, an injective endomorphism $\fr_n$ of $R$. One thus obtains a canonical action of the semigroup $\nt$ on any such $R$ and it is thus natural to work in the topos $\wnt$ of sets endowed with an action of $\nt$. 
\begin{defn}\label{site} The arithmetic site $\aarith=\arith$ is the topos $\wnt$
endowed with the  {\em structure sheaf} $\cO:=\zmax$ viewed as a semiring in the topos using the action of $\nt$ by the Frobenius endomorphisms.
\end{defn}
 The topological space underlying the arithmetic site is the Grothendieck topos of sets endowed with an action of  the multiplicative mono\"{\i}d $\nt$ of non-zero positive integers. 
 As we have seen above the semifield $\rmax$ of tropical real numbers   admits a one parameter group of  Frobenius
automorphisms $\fr_\lambda$, $\lambda\in \R_+^\times$, given by $\fr_\lambda(x)=x^\lambda$ $\forall x\in \rmax$. 
Using a straightforward extension in the context of semi-ringed topos of the classical notion of algebraic geometry of a point over a ring, one then gets the following result which gives the bridge between the noncommutative geometry and topos points of view:
\begin{thm}\cite{CC,CCas1} \label{structure3z} The set of points of the arithmetic site  $\aarith$ over $\rmax$ is canonically isomorphic with $X_\Q=\Q^\times\backslash \A_\Q/\hatz$. The action of the Frobenius
automorphisms $\fr_\lambda$ of $\rmax$ on these points corresponds to the action of the idele class group on $X_\Q=\Q^\times\backslash \A_\Q/\hatz$.
\end{thm}
The square of the arithmetic site is the topos $\wntb$ endowed with the structure sheaf defined globally by the multiplicatively cancellative semiring associated to the tensor square $\nbo$ over the  smallest Boolean semifield of characteristic one.  In this way one obtains the semiring whose elements are Newton polygons and whose operations are given by the convex hull of the union and the sum. The points of the  square of the arithmetic site over $\rmax$ coincide with the product of the points of the  arithmetic site over $\rmax$. Then, we describe the Frobenius correspondences $\Psi(\lambda)$ as congruences on the square parametrized by positive real numbers $\lambda\in \R_+^\times$.
 
\vsp
\centerline{\colorbox{lightgray}{\parbox[top][2.5cm][c]{13cm}{
The remarkable fact at this point is that while the arithmetic site is constructed as a combinatorial object of countable nature it  possesses nonetheless a one parameter semigroup of ``correspondences" which can be viewed as congruences in the square of the site.
}}}

\vsp
 In the context of semirings, the congruences \ie the equivalence relations compatible with addition and product, play the role of the ideals in ring theory. The Frobenius correspondences $\Psi(\lambda)$, for a rational value of $\lambda$, are deduced from the diagonal of the square, which is described by the product structure of the semiring, by  composition with the Frobenius endomorphisms. We interpret these correspondences geometrically, in terms of the congruence relation on Newton polygons corresponding to their belonging to the same half planes with rational slope $\lambda$. These congruences continue to make sense also for irrational values of $\lambda$ and are described using the best rational approximations of $\lambda$, while different values  of the parameter give rise to distinct congruences. 
 The composition of the Frobenius correspondences is  given for $\lambda, \lambda' \in \R_+^\times$  such that $\lambda\lambda'\notin \Q$ by the rule \cite{CC,CCas1}
\begin{equation}\label{compideps}
\Psi(\lambda)\circ \Psi(\lambda')=\Psi(\lambda\lambda')
\end{equation}
The same equality still holds if $\lambda$ and $\lambda'$ are rational numbers. When  $\lambda, \lambda'$ are irrational and $\lambda\lambda'\in \Q$ one has
\begin{equation}\label{compideps1}
\Psi(\lambda)\circ \Psi(\lambda')=\id_\epsilon\circ \Psi(\lambda\lambda')
\end{equation}
where $\id_\epsilon$ is the tangential deformation of the identity correspondence.

\subsubsection{The scaling site and Riemann-Roch theorems}
The Scaling Site $\scal1$, \cite{CCss},  is the algebraic geometric space  obtained from the arithmetic site $\aarith$ of \cite{CC,CCas1} by extension of scalars from the Boolean semifield $\B$ to the tropical semifield $\rmax$. The points of $\scal1$ are the same as the points $\aarith(\rmax)$ of the arithmetic site over $\rmax$. But $\scal1$ inherits from its structural sheaf a natural structure of tropical curve, in a generalized sense,  allowing one to define the sheaf of rational functions and to investigate an adequate version of the Riemann-Roch theorem in characteristic $1$. In \cite{CCss}, we tested this structure by restricting it to the periodic orbits of the scaling flow, \ie the points over the image of $\Spec\Z$ under the canonical morphism of toposes $\Theta:\spz\to \aarith$ (\cf~\cite{CCas1}, \S 5.1). We found that for each prime $p$ the corresponding circle of length $\log p$ is endowed with a quasi-tropical structure which turns this orbit into the analogue $C_p=\R_+^*/p^\Z$ of a classical elliptic curve $\C^*/q^\Z$. In particular  rational functions, divisors, etc all make sense. A new feature is that the degree of a divisor can now be any real number. The Jacobian of $C_p$ (\ie the quotient $J(C_p)$  of the group of divisors of degree $0$ by principal divisors)  is a cyclic group of order $p-1$. For each divisor $D$ there is a corresponding Riemann-Roch problem with solution space $H^0(D)$ and  the continuous dimension $\cdim(H^0(D))$ of this $\rma$-module is defined as the limit
\begin{equation}\label{rr1}
\cdim(H^0(D)):=\lim_{n\to \infty} p^{-n}\tdim(H^0(D)^{p^n})
\end{equation}
where $H^0(D)^{p^n}$ is a natural filtration and $\tdim(\cE)$ is  the topological dimension of an $\rma$-module $\cE$.  One has the following Riemann-Roch formula \cite{CCss},
\vsp
\begin{thm}\label{RRperiodic}
$(i)$~Let $D\in \div(C_p)$ be a divisor with $\deg(D)\geq 0$. Then the limit in \eqref{rr1} converges and one has  
$\cdim(H^0(D))=\deg(D)$.\newline
$(ii)$~The following Riemann-Roch formula holds
\begin{equation*}
\cdim(H^0(D))-\cdim(H^0(-D))=\deg(D)\qqq  D\in \div(C_p)
\end{equation*}
\end{thm}
The appearance of arbitrary positive real numbers as continuous dimensions in the  Riemann-Roch formula is due to the density in $\R$ of the subgroup $H_p\subset \Q$ of fractions with denominators a power of $p$.
 This outcome  is the analogue in characteristic $1$ of what happens for modules over matroid  $C^*$-algebras and the type II normalized dimensions  as in \cite{dix}. 

\vsp
\centerline{\colorbox{lightgray}{\parbox[top][3cm][c]{13cm}{
At this point, what is missing is an intersection theory and a  Riemann-Roch theorem on the square of the arithmetic site.
One expects that the  right hand side of the Riemann-Roch formula will be of the form $\frac 12 D.D=\inter(f,f)$ when  the divisor $D$ is of the form 
$$
D(f)=\int \Psi(\lambda) f(\lambda) d^*\lambda
$$
}}}
\vsp

Here $f(\lambda)$ is a real valued function with compact support of the variable $\lambda\in \R_+^*$ and $\inter(f,f)$ is as in \eqref{negcrit}. More precisely $D.D$ should be obtained as the intersection number of $D\circ \tilde D$ (defined using composition of correspondences) with the diagonal $\Delta$ and hence as a suitably defined distributional trace as for the counting function $N(u)$ of \S \ref{counting} so that $\frac 12 D(f).D(f)=\inter(f,f)$ with the notations of \eqref{negcrit}. So far the Riemann-Roch formula in tropical geometry is limited to curves and there is no Serre duality or good cohomological version of $H^j$ for $j\neq 0$, but in the above context one can hope that a Riemann-Roch inequality of the type \eqref{rrine}, \ie of the form
\begin{equation*}
\cdim(H^0(D))+\cdim(H^0(-D))\geq \frac 12 D.D
\end{equation*}
would suffice to apply the strategy of   Section \ref{rrstrat} to prove the key inequality \eqref{negcrit}.

\begin{table}[H]
\begin{center}
\caption{Here are a few entries in the analogy:}
\label{tab:1}       
\vsp\vsp
\begin{tabular}{| c | c |}
\hline
&\\
$C$ curve over $\F_q$ & Arithmetic Site $\aarith= ( \wnt,\zmax)$ over $\B$ \\ &\\
\hline &\\
Structure sheaf $\cO_{C}$& Structure sheaf $\zmax$  \\
&\\
\hline &\\
$\bar C=C\otimes_{\F_q} \bar\F_q$& Scaling Site $\scal1=(\rnt,\cO)$ over $\rmax$\\ 
&\\
\hline &\\
$C( \bar\F_q)=\bar C(\bar\F_q)$& $\aarith(\rmax)=\scal1(\rmax)$\\ &\\
\hline &\\
Galois action on $C( \bar\F_q)$& Galois action on $\aarith(\rmax)$ \\ &\\
\hline &\\
Structure sheaf $\cO_{\bar C}$& Structure sheaf $\cO=\zmax\hat\otimes_\B\rmax$\\
of $\bar C=C\otimes_{\F_q} \bar\F_q$ & piecewise affine convex functions, integral slopes\\ &\\
\hline&\\
Sheaf $\cK$ of rational functions  $\bar C$ & Sheaf $\cK$ of  piecewise affine functions\\ 
on $\bar C=C\otimes_{\F_q} \bar\F_q$& with integral slopes\\& \\
\hline&\\
Cartier divisors $=$ sections of $\cK/\cO^*$ & Sections of  $\cK/\cO^*$\\
&\\
\hline &\\
$X= \bar C\times \bar C$ & $\scal1\times \scal1$\\ &\\
\hline &\\
$D=\sum a_k\Psi^k$ &  $D=\int \Psi(\lambda) f(\lambda) d^*\lambda$ \\ 
Frobenius correspondence $\Psi$& Correspondences $\Psi(\lambda)$\\ &\\
\hline
\end{tabular}
\end{center}
\end{table}

\section{Absolute Algebra and the sphere spectrum}

Even if the Riemann-Roch strategy of Section \ref{algeomattack} happened to be successful, one should not view the arithmetic and scaling sites for more than what they are, namely a semiclassical shadow of a still mysterious structure dealing with compactifications of  $\Spec \Z$. An essential role in the unveiling of this structure should be played, for the reasons briefly explained below, by the discovery made by algebraic topologists in the 80's (see \cite{DGM}) that in their world of ``spectra" (in their sense) the sphere spectrum is a generalized ring $\sss$ which is more fundamental than the ring $\Z$ of integers, while 
the latter becomes an $\sss$-algebra. Over the years the technical complications of dealing with spaces ``up to homotopy" have greatly been simplified, in particular for the smash product of spectra. For the purpose of arithmetic applications, Segal's $\Gamma$-rings provide a very simple algebraic framework which succeeds to unify several attempts pursued in recent times in order to define the meaning of  ``absolute algebra".  In particular it contains the following three possible categories that had been considered previously to handle this unification: namely the category $\mathbf\cM$ of mono\"{\i}ds as in  \cite{deit, deit1,CC1,CC3}, the category $\mathbf\cH$ of hyperrings of \cite{CC2,CC4,CC5} and finally the category $\mathbf\cS$ of  semirings as in \cite{C,CC,CCas1,CCss}. Thanks to the work of L. Hesselholt and I. Madsen briefly explained below in \S \ref{topcyc} one now has at disposal a candidate cohomology theory in the arithmetic context: topological cyclic homology.

\subsection{Segal's $\Gamma$-rings}

  Let $\Gamma^{\rm op}$ be  the small, full subcategory  of the category of finite pointed sets whose objects are the the pointed finite sets\footnote{where $0$ is the base point.} $k_+:=\{0,\ldots ,k\}$, for $k\geq 0$.  The object $0_+$ is both initial and final so that $\gop$ is a  {\em pointed category}. The notion  of a discrete $\Gamma$-space, \ie of a $\Gamma$-set  is as follows:
\begin{defn}\label{defngamset} A $\Gamma$-set $F$ is a  functor $F:\gop\longrightarrow\Ses$ between pointed categories from $\gop$ to the category of pointed sets. \end{defn}
The morphisms $\Hom_\gop(M,N)$ between two $\Gamma$-sets are natural transformations of functors. 
The category $\gam$ of $\Gamma$-sets is a symmetric closed monoidal category (\cf~\cite{DGM}, Chapter II).  The monoidal structure is given by the smash product (denoted $X\wedge Y$) of $\Gamma$-sets which is a Day product. The closed structure property is shown in \cite{Lyd} (\cf~also \cite{DGM} Theorem 2.1.2.4). The specialization of Definition 2.1.4.1. of \cite{DGM} to the case of $\Gamma$-sets yields the following
\begin{defn}\label{defnsalg} A $\Gamma$-ring $\mathcal A$  is a $\Gamma$-set $\mathcal A: \gop\longrightarrow\Ses$  endowed with an associative multiplication 
$\mu:\cA \wedge \cA\to \cA$ and a unit~ $1: \sss\to \cA$, where  
$
\sss:\gop\longrightarrow \Ses
$
is the inclusion  functor.
\end{defn}
 Thus $\Gamma$-rings\footnote{equivalently $\sss$-algebras} make sense and the sphere spectrum corresponds to the simplest possible $\Gamma$-ring:  $\sss$. One can then easily identify the category $\gam$ of $\Gamma$-sets with the category $\catmo(\sss)$ of $\sss$-modules. 
In \cite{Durov}, N. Durov developed a geometry over $\F_1$ intended for Arakelov theory applications by using  monads as  generalizations of classical rings. 
While in the context of  \cite{Durov} the tensor product $\Z\otimes_{\F_1} \Z$ produces an uninteresting output isomorphic to $\Z$,  we showed in \cite{CCsalg}  that the same tensor square, re-understood in the theory of $\sss$-algebras, provides a highly non-trivial object. 
The Arakelov compactification of $\Spec\Z$ is endowed naturally with a structure sheaf of $\sss$-algebras and each Arakelov divisor provides a natural sheaf of modules over the structure sheaf. This new structure of $\overline{\Spec\Z}$ over $\sss$ endorses a one parameter group of weakly invertible sheaves whose tensor product rules are the same as   the composition rules \eqref{compideps}, \eqref{compideps1} of the Frobenius correspondences over
the arithmetic site \cite{CC, CCas1}.    The category $\catmo(\sss)$ of $\sss$-modules is not an abelian category and thus the tools of homological algebra need to be replaced along the line of the Dold-Kan correspondence, which for an abelian category $\cA$ gives the correspondence between chain complexes in $\geq 0$ degrees and simplicial objects \ie objects of $\cA^\dop$. 

\vsp
\centerline{\colorbox{lightgray}{\parbox[top][2.8cm][c]{13cm}{
At this point one has the following simple  but very important observation that $\Gamma$-spaces should be viewed as simplicial objects in $\gam\equiv \catmo(\sss)$, so that homotopy theory should be considered as the homological algebra corresponding to the ``absolute algebra" taking place over $\sss$.
}}}
\vsp

We refer to Table \ref{tab:2} for a short dictionary.
 The category  of $\Gamma$-spaces is the central tool of  \cite{DGM}, while the relations between algebraic $K$-theory and topological cyclic homology is the main topic.

\begin{table}[H]
\begin{center}
\caption{Short dictionary homology--homotopy}
\label{tab:2}       
\vsp\vsp
\begin{tabular}{| c | c |}
\hline &   \\
$X\in Ch_{\geq 0}(\cA)$ &  $M\in \catmo(\sss)^\dop$ \\
&  \\
\hline &   \\
$H_q(X)$&  $\pi_q(M)$ 
 \\ &  \\
\hline &   \\
$H_q(f):H_q(X)\simeq H_q(Y)$ &  $\pi_q(f):\pi_q(M)\simeq \pi_q(N)$    \\&    \\
quasi-isomorphism&   weak equivalence \\
&    \\
\hline &   \\
$f_n:X_n\stackrel{\subset}{\to} Y_n$ &  cofibration   \\
 +  projective cokernel &   (stable) \\
&    \\
\hline  &   \\
$f_n:X_n\to Y_n$ &  stable   \\
surjective if $n>0$&   fibration \\
&    \\
\hline  \end{tabular}
\end{center}
\end{table}

\subsection{Topological cyclic homology}\label{topcyc}

As shown in \cite{CCsalg} the various attempts done in recent times to develop  ``absolute algebra" are all unified by means of  the well established concept of $\sss$-algebra, \ie of $\Gamma$-rings.  Moreover (\cff\cite{DGM}) this latter notion is at the root of  the theory of topological cyclic homology which can be understood as cyclic homology over the absolute base $\sss$, provided one uses the appropriate Quillen model category. In particular, topological cyclic homology is now available to understand the new structure of $\overline{\Spec\Z}$ using its structure sheaf and modules. The use of cyclic homology in the arithmetic context is backed up by the following two results:
\begin{itemize}
\item  At the archimedean places, and after the initial work of Deninger  \cite{Den1,Den2} to recast the archimedean local factors of arithmetic varieties \cite{Se3} as regularized determinants,  we showed in  \cite{CC6} that cyclic homology in fact gives the correct infinite dimensional (co)homological theory for arithmetic varieties. The key operator $\Theta$ in this context is  the generator of the $\lambda$-operations $\Lambda(k)$ \cite{Loday,weibel,Weibelcris} in cyclic theory.
 More precisely, the action $u^\Theta$ of the multiplicative group $\R_+^\times$ generated by $\Theta$ on cyclic homology, is uniquely determined by its restriction to the dense subgroup $\Q_+^\times\subset \R_+^\times$ where it is  given by the formula
\begin{equation}\label{actiontheta}
    k^\Theta|_{HC_n}=\Lambda(k)\,k^{-n} \qqq n\geq 0, \,  \ k\in \N^\times\subset \R_+^\times
\end{equation}
 Let $X$ be a smooth, projective variety of dimension $d$ over an algebraic number field $\K$ and let $\nu\vert\infty$ be an archimedean place of $\K$.
Then, the action of the operator $\Theta$ on the archimedean cyclic homology $\har$ (\cff\cite{CC6}) of $X_\nu$ satisfies 
\begin{equation}\label{dettheta0}
    \prod_{0\leq w \leq 2d} L_\nu(H^w(X),s)^{(-1)^{w}}=\frac{det_\infty(\frac{1}{2\pi}(s-\Theta)|_{\har_{\rm od}(X_\nu)})}{
    det_\infty(\frac{1}{2\pi}(s-\Theta)|_{\har_{\rm ev}(X_\nu)})}
\end{equation}
The left-hand side of \eqref{dettheta0} is the product of Serre's archimedean local factors of the complex $L$-function of $X$ (\cf\cite{Se3}). On the right-hand side, $det_\infty$ denotes the regularized determinant  and one sets $$\har_{\rm ev}(X_\nu)=\bigoplus_{n=2k\ge 0} \har_{n}(X_\nu), \  \har_{\rm od}(X_\nu)=\bigoplus_{n=2k+1\ge 1} \har_{n}(X_\nu)$$

\item L. Hesselholt and I. Madsen have shown  (\cf~\eg \cite{HM,H,H1})  that the de Rham-Witt complex, an essential ingredient of crystalline cohomology (\cff\cite{Berth, Illusie}),  arises 
naturally when one studies the topological cyclic homology of smooth algebras over a perfect field of finite characteristic. 
One of the remarkable features in their work is that the arithmetic ingredients such as the Frobenius and restriction maps are naturally present in the framework of topological cyclic homology. Moreover L. Hesselholt has shown \cite{Hessel1} how topological periodic cyclic homology with its inverse Frobenius operator may be used to give a cohomological interpretation of the Hasse-Weil zeta function of a scheme smooth and proper over a finite field in the form (\cff\cite{Hessel1}):
\begin{equation}\label{detthetaH}
    \zeta(X,s)=\frac{det_\infty(\frac{1}{2\pi}(s-\Theta)|_{\hes_{\rm od}(X)})}{
    det_\infty(\frac{1}{2\pi}(s-\Theta)|_{\hes_{\rm ev}(X)})}
\end{equation}
\end{itemize}

\vsp
\centerline{\colorbox{lightgray}{\parbox[top][2.9cm][c]{13cm}{
The similarity between \eqref{dettheta0}  and \eqref{detthetaH} (applied to a place of good reduction) suggests the existence of a global formula for the $L$-functions of arithmetic varieties, involving cyclic homology of $\sss$-algebras, and of a Lefschetz formula in which the local factors appear from the periodic orbits of the action of $\R_+^*$. 
}}}
\vsp

One of the stumbling blocks in order to reach a satisfactory cohomology theory  is the problem of coefficients. Indeed, the natural coefficients at a prime $p$ for crystalline cohomology are an extension of $\Q_p$ and it is traditional to relate them with complex numbers by an embedding of fields. Similarly, \eqref{detthetaH} uses an embedding of the Witt ring $W(\F_q)\to \C$. To an analyst it is clear that since such embeddings cannot be measurable\footnote{A measurable group homomorphism from $\Z_p^\times$ to $\C^\times$ cannot be injective} they will never be effectively constructed. This begs for a better construction, along the lines of Quillen's computation of the algebraic $K$-theory of finite fields, which instead would only involve the ingredient of the Brauer lifting, \ie a group injection of the multiplicative group  of $\bar \F_p$ as roots of unity in $\C$.

\subsection{Final remarks}
The Riemann hypothesis has been extended far beyond its original formulation to the question of localization of the zeros of $L$-functions. There are  a number of constructions of $L$-functions coming from three different sources, Galois representations, automorphic forms and arithmetic varieties. Andr\' e Weil liked to compare (\cff \cite{B2} \S 12 and also \cite{weilcomplete} vol. 1, p. 244--255 and vol. 2, p. 408--412), the puzzle of these three different writings to the task of deciphering hieroglyphics with the help of the Rosetta Stone. In some sense the $L$-functions  play a role in modern mathematics similar to the role of polynomials in ancient mathematics, while the explicit formulas play the role of the expression of the symmetric functions of the roots in terms of the coefficients of the polynomial. If one follows this line of thought, the RH should be seen only as a first step since in the case of polynomials  there  is no way one should feel to have understood the zeros  
once one proves that they  are, say, real numbers. In fact Galois formulated precisely the problem as that of finding all numerical relations between the roots of an equation, with the  trivial ones being given by the symmetric functions, while the others, when determined, will reveal 
a complete understanding of the zeros as obtained, in the case of polynomials, by Galois theory.
In a fragment, page 103, of the complete works of Galois  \cite{Galois} concerning the memoir of February 1830, he delivers the essence of his theory:
\begin{quote}{\it 
Remarquons que tout ce qu'une \'equation num\'erique peut avoir de particulier, doit provenir de certaines relations entre les racines. Ces relations seront rationnelles c'est-\`a-dire qu'elles ne contiendront d'irrationnelles que les coefficients de l'\'equation et les quantit\'es adjointes. De plus ces relations ne devront pas \^etre invariables par toute substitution op\'er\'ee sur les racines, sans quoi on n'aurait rien de plus que dans les \'equations litt\'erales. Ce qu'il importe donc de conna\^itre, c'est par quelles substitutions peuvent \^etre invariables des relations entre les racines, ou ce qui revient au m\^eme, des fonctions des racines dont la valeur num\'erique est d\'eterminable rationnellement.\footnote{In 2012 I had to give, in the French academy of Sciences, the talk devoted to the 200-th anniversary of the birth of Evariste Galois. On that occasion I read for the $n+1$-th time the book of his collected works and was struck by the pertinence of the above quote in the analogy with $L$-functions. In the case of function fields one is dealing with Weil numbers and one knows a lot on their Galois theory using results such as those of Honda and Tate \cff\cite{Tatebbk}.}  
}\end{quote}

\subsubsection*{Acknowledgement} I am grateful to J. B. Bost for  the reference \cite{Tatebbk}, to J. B. Bost, P. Cartier, C. Consani, D. Goss, H. Moscovici, M. Th. Rassias, C. Skau and W. van Suijlekom for their detailed comments, to S. Gaubert for his help in Section 4.2 and to Lars Hesselholt for his comments and for allowing me to mention his forthcoming paper \cite{Hessel1}.

\end{document}